\documentclass[12pt,a4paper,final]{amsart}
\usepackage{version}
\usepackage[notcite,notref]{showkeys}
\usepackage{amsmath,amsthm,amssymb,latexsym,epsfig,graphicx,subfigure}
\numberwithin{equation}{section}

\newtheorem{thm}{Theorem}[section]

\theoremstyle{definition}

\theoremstyle{definition}

\newcommand{\Rn}{\mathbb{R}^{n}}

\newcommand{\R}{\mathbb{R}}

\renewcommand{\H}{\mathbb{H}}
\renewcommand{\H}{\mathbb{H}}

\newcommand {\grtrsim} {\ {\raise-.5ex\hbox{$\buildrel>\over\sim$}}\ }
\newcommand{\e}{\varepsilon}

\newcommand{\khii}{\text{\lower -.4ex\hbox{$\chi$}}}
\DeclareMathOperator{\spt}{spt}

\renewcommand{\a}{\alpha}

\newcommand{\restrict}{\begin{picture}(12,12)
                       \put(2,0){\line(1,0){8}}
                       \put(2,0){\line(0,1){8}}
                      \end{picture}}

\begin{document}
\title {A survey on the Hausdorff dimension of intersections}
\author{Pertti Mattila}

\thanks{} 

\begin{abstract} 
Let $A$ and $B$ be Borel subsets of the Euclidean $n$-space with $\dim A + \dim B > n$. This is a survey on the question: what can we say about the Hausdorff dimension of the intersections $A\cap (g(B)+z)$ for generic orthogonal transformations $g$ and translations by $z$. 
\end{abstract}

\maketitle

In this revision Theorems \ref{level} and \ref{level1} are modified to Theorems \ref{levela} and \ref{level1a} by replacing the assumptions \eqref{L2} and 
\eqref{L3} by \eqref{L2a} and \eqref{L3a} and hence by removing the positive lower density assumptions from Theorems \ref{level} and \ref{level1}. As a consequence, 
the positive lower density assumption is also removed from Theorem \ref{inter}. The possibility of these improvements was anticipated by Harris in \cite{H2}. I am grateful to Terence Harris for useful comments.

\section{Introduction} \footnote{This survey is based on the talk I gave in Karoly Simon's 60+1 birthday conference in Budapest in June 2022.}
The books \cite{M4} and \cite{M5} contain most of the required background information and the proofs of some of the results discussed below.

Let $\mathcal L^n$ stand for the Lebesgue measure on the Euclidean $n$-space $\Rn$ and let $\dim$ stand for the Hausdorff dimension and $\mathcal H^s$ for $s$-dimensional Hausdorff measure. For $A\subset \Rn$, denote by $\mathcal M(A)$ the set of Borel measures $\mu$ with $0<\mu(A)<\infty$ and with the compact support $\spt\mu\subset A$.

We let $O(n)$ denote the orthogonal group of $\Rn$ and $\theta_n$ its Haar probability measure.  The main fact needed about the measure $\theta_n$ is the inequality: 

\begin{equation}\label{theta}
\theta_n(\{g\in O(n):|x-g(z)|<r\})\lesssim (r/|z|)^{n-1}\ \text{for}\ x,z\in\Rn, r>0.
\end{equation}
This is quite easy, in fact trivial in the plane.

Let $A$ and $B$ be Borel subsets of $\Rn$ with Hausdorff dimensions $s=\dim A$ and $t=\dim B$. What can we say about the Hausdorff dimensions of the intersections of $A$ and typical rigid motions of $B$? More precisely, of $\dim A\cap (g(B)+z)$ for almost all $g\in O(n)$ and for $z\in\Rn$ in a set of positive Lebesgue measure. Optimally one could hope that this dimension is given by the bigger of the numbers $s+t-n$ and 0, which happens when smooth surfaces meet in a general position.

The problem on the upper bound is much easier than on the lower bound. Let 
\begin{equation}\label{eq0}
V_z = \{(x,y)\in\Rn\times\Rn: x=y+z\},\ z\in\Rn,\end{equation}
be the $z$ translate of the diagonal in $\Rn\times\Rn$, and let $\pi$ be the projection $\pi(x,y)=x$. Then 
\begin{equation}\label{eq1}
A\cap (g(B)+z) = \pi((A\times g(B))\cap V_z),
\end{equation}
and it follows from a Fubini-type inequality for Hausdorff dimension, \cite[Theorem 7.7]{M4}, that for any $g\in O(n)$,
\begin{equation}\label{eq2}
\dim A\cap (g(B)+z) \leq \dim(A\times B) - n\ \text{for almost all}\ z\in\Rn,
\end{equation}
provided $\dim(A\times B) \geq n$. We have always $\dim(A\times B) \geq \dim A +\dim B$ and the equation $\dim(A\times B) = \dim A +\dim B$ holds if, for example, $0<\mathcal H^s(A)<\infty, 0<\mathcal H^t(B)<\infty$, and one of the sets has positive lower density, say
\begin{equation}\label{eq3}
\theta_{\ast}^s(A,x)=\liminf_{r\to 0}r^{-s}\mathcal H^s(A\cap B(x,r))>0\ \text{for}\ \mathcal H^s\ \text{almost all}\ x\in A. 
\end{equation}
Even the weaker condition that the Hausdorff and packing dimensions of $A$ agree suffices, see \cite{M4}, pp. 115-116. Then we have
\begin{equation}\label{eq4}
\dim A\cap (g(B)+z) \leq \dim A +\dim B - n\ \text{for almost all}\ z\in\Rn,
\end{equation}
provided $\dim A +\dim B \geq n$. Without some extra condition this inequality fails badly: for any $0\leq s\leq n$ there exists a Borel set $A\subset \Rn$ of dimension $s$ such that 
$\dim A\cap f(A)=s$ for all similarity maps $f$ of $\Rn$. This was proved by Falconer in \cite{F3}, see also Example 13.19 in \cite{M4} and the further references given there.  
 
We have the lower bound for the dimension of intersections if we use larger transformation groups, for example similarities:

\begin{thm}\label{thm1}
Let $A$ and $B$ be Borel subsets of $\Rn$ with $\dim A + \dim B > n$. Then for every $\e>0$,
$$\mathcal L^n(\{z\in\Rn:\dim A\cap (rg(B)+z) \geq \dim A + \dim B - n - \e\}) > 0,$$
for almost all $g\in O(n)$ and almost all $r > 0$.
\end{thm}

If $A$ and $B$ have positive and finite Hausdorff measure, $\e$ is not needed. This theorem was proved in the 1980s independently by Kahane \cite{K} and in \cite{M1}. More generally, Kahane proved that the similarities can be replaced by any  closed subgroup of the general linear group of $\Rn$ which is transitive outside the origin. He gave applications to multiple points of stochastic processes.

There are many special cases where the equality $\dim A\cap (g(B)+z) = \dim A + \dim B - n$ holds for almost all $g$ and for $z$ in a set of positive measure. The case where one of the sets is a plane,  initiated by Marstrand in \cite{M}, has been studied a lot, see discussions in \cite[Chapter 10]{M4} and \cite[Chapter 6]{M5}, and \cite{MO} for a more recent result. More generally, one of the sets can be rectifiable, see \cite{M1}.

The main open problem is: what conditions on the Hausdorff dimensions or measures of $A$ and $B$ guarantee that for $\theta_n$ almost all $g\in O(n)$,

\begin{equation}\label{eq5}
\mathcal L^n(\{z\in\Rn:\dim A\cap (g(B)+z) \geq \dim A + \dim B - n\}) > 0,
\end{equation}
or perhaps for all $\e>0$,
\begin{equation}\label{eq6}
\mathcal L^n(\{z\in\Rn:\dim A\cap (g(B)+z) \geq \dim A + \dim B - n-\e\}) > 0?
\end{equation}
If one of the sets is a Salem set, that is, it supports a measure with an optimal Fourier decay allowed by its Hausdorff dimension, then \eqref{eq6} holds without dimensional restrictions, see \cite{M3}. I expect \eqref{eq6} to be true for all Borel subsets $A$ and $B$ of $\Rn$.

Below I shall discuss some partial results on this question. I shall also say something about the exceptional sets of transformations. 
%But first I shall discuss the methods.

In this survey I shall concentrate on Hausdorff dimension and general Borel sets. For remarks and references about related results on other dimensions, see \cite[Section 13.20]{M4} and \cite[Section 7.3]{M5}. There is a rich literature on various questions about intersections of dynamically generated and related sets. For recent results and further references, see \cite{S1}, \cite{Wu}, \cite{Y}. For probabilistic sets, see \cite{SS} and its references.

I would like to thank the referees for useful comments.

\section{Projections and plane intersections}\label{Projections}

This topic can be  thought of as a study of integral-geometric properties of fractal sets and Hausdorff dimension.  Let us briefly review some of the basic related results on projections and plane sections. This  was started by Marstrand in \cite{M} in the plane. His main results in general dimensions are the following. Let $G(n,m)$ be the Grassmannian of linear $m$-dimensional subspaces of $\Rn$ and $P_V:\Rn\to V$ the orthogonal projection onto $V\in G(n,m)$. Let also $\gamma_{n,m}$ be the orthogonally invariant Borel probability measure on $G(n,m)$.

\begin{thm}\label{mar1}
Let $A\subset\Rn$ be a Borel set. Then for almost all $V\in G(n,m)$,
\begin{equation}\label{mareq1}
\dim P_V(A) = \dim A\ \text{if}\ \dim A \leq m,
\end{equation}
and 
\begin{equation}\label{mareq12}
\mathcal H^m(P_V(A))>0\ \text{if}\ \dim A > m.
\end{equation}
\end{thm}

\begin{thm}\label{mar2}
Let $n-m \leq s \leq n$ and let $A\subset\Rn$ be $\mathcal H^s$ measurable with $0<\mathcal H^s(A)<\infty$. Then for almost all $V\in G(n,m)$,
\begin{equation}\label{mareq3}
\mathcal H^{n-m}(\{u\in V^{\perp}: \dim (A\cap (V+u)) = s+m-n\})>0,
\end{equation}
and for almost all $V\in G(n,n-m)$ and for $\mathcal H^s$ almost all $x\in A$,
\begin{equation}\label{mareq4}
\dim(A\cap (V+x)) = s+m-n.
\end{equation}
\end{thm}

One can sharpen these results by deriving estimates on the Hausdorff dimension of the exceptional sets of the planes $V$. For the first part of Theorem \ref{mar1} this was first done by Kaufman in \cite{Ka} in the plane, then in \cite{M0} and \cite{KM} in higher dimensions. For the second part of Theorem \ref{mar1} the exceptional set estimates were proven by Falconer in \cite{F1}. Thus we have, recall that $\dim G(n,m)= m(n-m)$:

\begin{thm}\label{exmar1}
Let $A\subset\Rn$ be a Borel set with $s=\dim A$. Then 
\begin{equation}\label{exmareq1}
\dim\{V\in G(n,m): P_V(A) < \dim A\}\leq s-m +m(n-m)\ \text{if}\ s\leq m,
\end{equation}
and 
\begin{equation}\label{exmareq2}
\dim\{V\in G(n,m): \mathcal H^m(P_V(A))=0\})\leq m-s+m(n-m)\ \text{if}\ s> m,
\end{equation}
\end{thm}

These inequalities are sharp by the examples in \cite{KM} (and their modifications), but the proof for \eqref{exmareq1} also gives the upper bound $t-m +m(n-m)$ if $\dim A$ on the left hand side is replaced by $t, 0\leq t \leq \dim A$. Then for $t <\dim A$ this is not always sharp, see the discussion in \cite[Section 5.4]{M5}.

For the plane sections Orponen proved in \cite{Or1}, see also \cite[Theorem 6.7]{M5}, the exceptional set estimate (which of course is sharp, as \eqref{exmareq2} is):

\begin{thm}\label{exmar2}
Let $n-m \leq s \leq n$ and let $A\subset\Rn$ be $\mathcal H^s$ measurable with $0<\mathcal H^s(A)<\infty$. Then there is a Borel set $E\subset G(n,m)$ such that 
$$\dim E\leq n-m-s+m(n-m)$$ 
and for $V\in G(n,m)\setminus E,$
\begin{equation}\label{mareq3}
\mathcal H^{n-m}(\{u\in V^{\perp}: \dim (A\cap (V+u)) = s+m-n\})>0.
\end{equation}
\end{thm}

We can also ask for exceptional set estimates corresponding to \eqref{mareq4}. We proved with Orponen \cite{MO} the following:

\begin{thm}\label{exmar3}
Let $n-m \leq s \leq n$ and let $A\subset\Rn$ be $\mathcal H^s$ measurable with $0<\mathcal H^s(A)<\infty$. Then the set $B$ of points $x\in\Rn$ with
$$\gamma_{n,m}(\{V\in G(n,m):\dim A\cap (V+x)=s+m-n\})=0$$
has dimension $\dim B\leq n-m$.
\end{thm}

Very likely, the bound $n-m$ is not sharp. When $m=1$, probably the sharp bound should be $2(n-1) - s$ in accordance with Orponen's sharp result for radial projections in \cite{Or2}.

Another open question is whether there could be some sort of non-trivial estimate for the dimension of the exceptional pairs $(x,V)$. 

\section{Some words about the methods}\label{methods}

The methods in all cases use Frostman measures. Suppose that the Hausdorff measures $\mathcal H^s(A)$ and $\mathcal H^t(B)$ are positive. Then there are 
$\mu\in\mathcal M(A)$ and $\nu\in\mathcal M(B)$ such that $\mu(B(x,r)) \leq r^s$ and $\nu(B(x,r)) \leq r^t$ for $x\in\Rn, r>0$. In particular, for 
$0<s<\dim A$ and $0<t<\dim B$ there are $\mu\in\mathcal M(A)$ and $\nu\in\mathcal M(B)$ such that $I_s(\mu)<\infty$ and $I_t(\nu)<\infty$, where the  $s$ energy $I_s(\mu)$ is defined by
$$I_s(\mu)=\iint |x-y|^{-s}\,d\mu x\,d\mu y.$$
Then the goal is to find intersection measures $\lambda_{g,z}\in\mathcal M(A\cap (g(B)+z))$ such that 
%\begin{itemize}
\begin{equation}\label{eq7} \spt\lambda_{g,z}\subset\spt\mu\cap(g(\spt\nu)+z),\end{equation}
\begin{equation}\label{eq8}\int\lambda_{g,z}(\Rn)\,d\mathcal L^nz = \mu(\Rn)\nu(\Rn)\ \text{for}\ \theta_n\ \text{almost all}\ g\in O(n),\end{equation}
\begin{equation}\label{eq9}\iint I_{s+t-n}(\lambda_{g,z})\,d\mathcal L^nz\,d\theta_n g\lesssim I_s(\mu)I_t(\nu).\end{equation} 
%\end{itemize} 
This would give \eqref{eq6}.

There are two closely related methods to produce these measures. The first, used in \cite{M1}, is based on \eqref{eq1}:  the intersections $A\cap (g(B)+z)$ can be realized as level sets of the projections $S_g$:
\begin{equation}\label{eq23}
S_g(x,y)=x-g(y),\ x,y\in\Rn,
\end{equation}
\begin{equation}\label{eq16.2}
A\cap (g(B)+z) = \pi ((A\times g(B))\cap S_g^{-1}\{z\}),\ \pi(x,y)=x.
\end{equation}
Notice that the map $S_g$ is essentially the orthogonal projection onto the $n$-plane $\{(x,-g(x)):x\in\Rn\}$.

Thus one slices (disintegrates) $\mu\times g_{\#}\nu$ ($g_{\#}\nu$ is the push-forward) with the planes $V_z=\{(x,y): x=y+z\}, z\in\Rn.$
For this to work, one needs to know that 
\begin{equation}\label{eq10}S_{g\#}(\mu\times\nu)\ll\mathcal L^n\ \text{for}\ \theta_n\ \text{almost all}\ g\in O(n).\end{equation}
This is usually proved by establishing the $L^2$ estimate 
\begin{equation}\label{eq16}
\iint S_{g_\#}(\mu\times\nu)(x)^2\,dx\,d\theta_ng \lesssim 1,
\end{equation}
which, by Plancherel's formula, is equivalent to
\begin{equation}\label{eq16.1}
\iint \mathcal F (S_{g\#}(\mu\times\nu))(x)^2\,dx\,d\theta_ng \lesssim 1,
\end{equation}
where $\mathcal F$ stands for the Fourier transform. 

The second method, used in \cite{K}, is based on convolution approximation. Letting $\psi_{\e}, \e>0$, be a standard approximate identity, set
$\nu_{\e}=\psi_{\e}\ast\nu$ and

\begin{equation}\label{eq11}
\nu_{g,z,\e}(x)=\nu_{\e}(g^{-1}(x-z)),\ x\in\Rn. 
\end{equation}
Then the plan is to show that when $\e\to 0$, the measures $\nu_{g,z,\e}\mu$ converge weakly to the desired intersection measures.

No Fourier transform is needed to prove Theorem \ref{thm1}, but the proofs of all theorems discussed below, except Theorems \ref{level}, \ref{levela}, \ref{level1} and \ref{level1a}, rely on the Fourier transform defined by
$$\widehat{\mu}(x)=\int e^{-2\pi ix\cdot y}\,d\mu y,~ x\in\Rn.$$
The basic reason for its usefulness in this connection is the formula
\begin{equation}\label{eq12}
I_s(\mu)=\iint|x-y|^{-s}\,d\mu x\,d\mu y=c(n,s)\int|\widehat{\mu}(x)|^2|x|^{s-n}\,dx,
\end{equation} 
which is a consequence of Parseval's formula and the fact that the distributional Fourier transform of the Riesz kernel $k_s, k_s(x)=|x|^{-s}$, is a constant multiple of $k_{n-s}$. 

The decay of the spherical averages, 
$$\sigma(\mu)(r)=r^{1-n}\int_{S(r)}|\widehat{\mu}(x)|^2\,d\sigma_r^{n-1}x, r>0,$$
 of $\mu\in \mathcal M(\Rn)$, where $\sigma_r^{n-1}$ is the surface measure on the sphere $S(r)=\{x\in\Rn:|x|=r\}$, often plays an important role.  
By integration in polar coordinates, if $\sigma(\mu)(r)\lesssim r^{-t}$ for $r>0$ and for some $t>0$, then $I_s(\mu)<\infty$ for $0<s<t$. Hence the best decay we can hope for under the finite $s$ energy assumption (or the Frostman assumption $\mu(B(x,r))\leq r^s)$) is $r^{-s}$. This is true when $s\leq (n-1)/2,$ see \cite[Lemma 3.5]{M5}, but false otherwise.

The decay estimates for $\sigma(\mu)(r)$ have been studied by many people, discussion can be found in \cite[Chapter 15]{M5}. The  best known estimates, due to Wolff, \cite{W}, when $n=2$ (the proof can also be found in \cite[Chapter 16]{M5}) and to Du and Zhang, \cite{DZ}, in the general case, are the following: Let  $\mu\in\mathcal M(\Rn), \spt\mu\subset B(0,1)$, with $\mu(B(x,r))\leq r^s$ for $x\in\R^n, r>0$. Then
for all $\e > 0, r>1$, 

\begin{equation}\label{dz}
\sigma(\mu)(r) \lesssim \begin{cases} r^{-(n-1)s/n+\e}\mu(B(0,1))\ \text{for all}\ 0<s<n,\\
r^{-(n-1)/2+\e}\mu(B(0,1))\ \text{if}\ (n-1)/2\leq s\leq n/2,\\
 r^{-s+\e}\mu(B(0,1))\ \text{if}\ 0<s\leq (n-1)/2.\\
\end{cases}
\end{equation}
The factor $\mu(B(0,1))$ can be checked from the proofs in \cite{W}, \cite{DZ} and \cite{M5}.
%for $n=2$ is explicitly stated in \cite[Theorem 16.1]{M5}. For $n>2$ it can be checked from the proof in \cite{DZ}.

The essential case for the first estimate is $s>n/2$, otherwise the second and third are better. Up to $\e$ these estimates are sharp when $n=2$. When $n\geq 3$ the sharp bounds are not known for all $s$, see \cite{D} for discussion and the most recent examples. As mentioned above, the last bound is always sharp.

\section{The first theorem}

If one of the sets has dimension bigger than $(n+1)/2$ we have the following theorem. It was proved in \cite{M2}, see also \cite[Theorem 13.11]{M4} or \cite[Theorem 7.4]{M5}: 

\begin{thm}\label{thm2}
Let $s$ and $t$ be positive numbers with $s+t > n$ and $s>(n+1)/2$. Let $A$ and $B$ be Borel subsets of $\Rn$ with $\mathcal H^s(A)>0$ and $\mathcal H^t(B)>0$. Then 
\begin{equation}\label{eq13}
\mathcal L^n(\{z\in\Rn:\dim A\cap (g(B)+z) \geq \dim A + \dim B - n\}) > 0
\end{equation}
for almost all $g\in O(n)$.
\end{thm}

The proof is based on the slicing method. The key estimate is  
\begin{equation}\label{eq14}
\mu\times\mu(\{(x,y): r\leq |x-y| \leq r+\delta\})\lesssim I_s(\mu)\delta r^{s-1}\end{equation}
if $\mu \in\mathcal M(\Rn), 0<\delta\leq r$ and $(n+1)/2\leq s<n$. This is combined with the inequality \eqref{theta}.

The inequality \eqref{eq14} is obtained with the help of the Fourier transform, and that is the only place in the proof of Theorem \ref{thm2} where the Fourier transform is needed.

One problem of extending Theorem \ref{thm2} below the dimension bound $(n+1)/2$ is that the estimate \eqref{eq14} then fails, at least in the plane by \cite[Example 4.9]{M5} and in $\R^3$ by \cite{IS}. 

In Section \ref{Excestimates} we discuss estimates on the exceptional sets of orthogonal transformations. The proof of Theorem \ref{extheo1} gives another proof for Theorem \ref{thm2} but under the stronger assumption $s+t>n+1$. On the other hand, Theorem \ref{inter} below holds with the assumption $s+(n-1)t/n>n$ but under the additional condition of positive lower density. Of course, $s+(n-1)t/n>n$ is sometimes stronger and sometimes weaker than $s>(n+1)/2, s+t>n$. For example, consider these in the plane. When $s=t$, the first one says $s>4/3$ and the second one $s>3/2$. On the other hand, when $s$ is slightly bigger than $3/2$, the first requires $t$ to be at least about 1, but the second allows $t=1/2$.  

Theorem \ref{thm2} says nothing in $\R^1$, and there is nothing to say: in \cite{M1} I constructed compact sets $A, B\subset\R$ such that $\dim A = \dim B = 1$ and $A\cap (B+z)$ contains at most one point for any $z\in\R$. With $A, B\subset\R$ as above, the $n$-fold Cartesian products $A\times\dots\times A\subset\Rn$ and $B\times\dots\times B\subset\Rn$ yield the corresponding examples in $\Rn$, that is, just with translations we get nothing in general.

Donoven and Falconer proved in \cite{DF} an analogue of Theorem \ref{thm2} for the isometries of the Cantor space. They didn't need any dimensional restrictions. They used martingales to construct the desired random measures with finite energy integrals on the intersections. 

\section{The projections $S_g$}

We now discuss a bit more the projections $S_g$, recall \eqref{eq23}. They are particular cases of restricted projections, which recently have been studied extensively, see \cite[Section 5.4]{M5}, \cite{M8}, \cite{H2} and \cite{GGGHMW} and the references given there. Restricted means that we are considering a lower dimensional subspace of the Grassmannian $G(2n,n)$. For the full Grassmannian we have Marstrand's projection theorem \ref{mar1}. 

As mentioned above, to prove Theorem \ref{thm2} one first needs to know \eqref{eq10} when $s+t > n$ and $s>(n+1)/2$ and $\mu$ and $\nu$ have finite $s$ and $t$ energies. A  simple proof using spherical averages is given in \cite[Lemma 7.1]{M5}. This immediately yields the weaker result: with the assumptions of Theorem \ref{thm2}, for almost all $g\in O(n)$,
\begin{equation}\label{eq15}
\mathcal L^n(\{z\in\Rn: A\cap (g(B)+z) \not= \emptyset\}) > 0,
\end{equation}
because \eqref{eq15} is equivalent to $\mathcal L^n(S_g(A\times B))>0$. 
Even for this I don't know if the assumption $s>(n+1)/2$ is needed. 

Let us first look at general Borel subsets of $\R^{2n}$:

\begin{thm}\label{thm3}
Let $A\subset\R^{2n}$  be a Borel set. If $\dim A > n+1$, then $\mathcal L^n(S_g(A))>0$ for $\theta_{n}$ almost all $g\in O(n)$. \end{thm}

This was proved in \cite{M8}. That paper also contains dimension estimates for $S_g(A)$ when $\dim A \leq n+1$ and estimates on the dimension of exceptional sets of transformations $g$. In particular, if $n\leq \dim A \leq n+1$, then
\begin{equation}\label{eq22}
\dim S_g(A)\geq \dim A-1\ \text{for}\ \theta_{n}\ \text{almost all}\ g\in O(n).\end{equation}

The bound $n+1$ in Theorem \ref{thm3} is sharp. This was shown by Harris in \cite{H1}. First, \eqref{eq22} is sharp.  The example for $\dim A=n$ is simply the diagonal $D=\{(x,x):x\in\Rn\}$. To see this suppose that $g\in O(n)$ is such that $\det g = (-1)^{n+1}$, which is satisfied by half of the orthogonal transformations. Then by some linear algebra $g$ has a fixed point, whence the kernel of $x\mapsto S_g(x,x)$ is non-trivial, so $\dim S_g(D)\leq n-1$. Taking the Cartesian product of $D$ with a one-dimensional set of zero $\mathcal H^1$ measure, we obtain $A$ with $\dim A = n+1$ and $\mathcal L^n(S_g(A))=0$, which proves the sharpness.

But this only gives an example $A$ of dimension $n+1$ for which $\mathcal L^n(S_g(A)) = 0$ for $g\in O(n)$ with measure $1/2$. Is there a counter-example that works for almost all $g\in O(n)$?
%On the other hand, if $\det g = (-1)^{n}$, there is a unit vector $e$ such that $g(e) = -e$, that is, $S_g(e,-e)=0$. Then, as above, we get from $E=\{(x,-x):x\in\Rn\}$ a closed set $B$ with $\dim B = n+1$ and 
%$\mathcal L^n(S_g(B))=0$. I believe (but I havn't checked) that a counter-example for all $g\in O(n)$ could be obtained performing a Cantor construction with cubes so that the set inside the first level cubes looks like $A$, inside the next level cubes like $B$, and so on.

Here are the basic ingredients of the proof of Theorem \ref{thm3}. They were inspired by Oberlin's paper \cite{O}. 

Let $0<n+1<s<\dim A$ and  $\mu\in\mathcal M(A)$ with $I_s(\mu) < \infty$, and let $\mu_g\in\mathcal M(S_g(A))$ be the push-forward of $\mu$ under $S_g$. The Fourier transform of $\mu_g$ is given by 
\begin{equation*}
\widehat{\mu_g}(\xi)=\widehat{\mu}(\xi,-g^{-1}(\xi)).
\end{equation*}
By fairly standard arguments, using also the inequality \eqref{theta}, one can then show that for $R>1$,
\begin{equation}\label{obeq}
\iint_{R\leq|\xi|\leq 2R}|\widehat{\mu}(\xi,-g^{-1}(\xi))|^2\,d\xi\,d\theta_n g\lesssim R^{n+1-s}.
\end{equation}
This is summed over the dyadic annuli, $R=2^k, k=1,2,\dots$. The sum converges since $s>n+1$. Hence for $\theta_n$ almost all $g\in O(n)$,\ $\mu_g$ is absolutely continuous with $L^2$ density, and so $\mathcal L^n(S_g(A))>0$. 

For product sets we can improve this, which is essential for the applications to intersections:
 
\begin{thm}\label{thm4}
Let $A, B\subset\R^{n}$  be Borel sets. If $\dim A + (n-1)\dim B/n > n$ or $\dim A+\dim B>n$ and $\dim A > (n+1)/2$, then $\mathcal{L}^n(S_g(A\times B))>0$ for $\theta_{n}$ almost all $g\in O(n)$.
\end{thm}

The case $\dim A > (n+1)/2$ is a special case of Theorem \ref{thm2}, recall \eqref{eq15}. The proof of the case $\dim A + (n-1)\dim B/n > n$ is based on the spherical averages and the first estimate of \eqref{dz}. Here is a sketch.

Suppose $A, B \subset B(0,1)$. Let $0<s<\dim A, 0<t<\dim B$ and $\e>0$ such that $s+(n-1)t/n-\e>n$, and let $\mu\in\mathcal M(A),\nu\in\mathcal M(B)$ with $\mu(B(x,r))\leq r^{s}, \nu(B(x,r))\leq r^{t}$ for $x\in\R^n, r>0$. Let $\lambda_g = S_{g\#}(\mu\times\nu)\in\mathcal M(S_g(A\times B)).$ Then $\widehat{\lambda_g}(\xi)=\widehat{\mu}(\xi)\widehat{\nu}(-g^{-1}(\xi))$. By \eqref{dz} we have
\begin{equation}\label{eq16b}
\begin{split}
&\iint|\widehat{\lambda_g}(\xi)|^2\,d\xi\,d\theta g
 =\int|\widehat{\mu}(\xi)|^2\sigma(\nu)(|\xi|)\,d\xi\\
&\lesssim\int|\widehat{\mu}(\xi)|^2|\xi|^{-(n-1)t/n+\e}\,d\xi\nu(B(0,1))\\
&=cI_{n-(n-1)t/n+ \e}(\mu)\nu(B(0,1))\lesssim \mu(B(0,1))\nu(B(0,1))<\infty.
\end{split}
\end{equation}
The last inequality is easy, see \cite[page 19]{M5}. This gives Theorem \ref{thm4}, and a quantitative estimate as in \eqref{L3a}.  

%In fact, for some results on the intersections below we again need absolute continuity as in \eqref{eq10}. In the case $s+(n-1)t>n$ we need the quantitative  estimate: if $s+(n-1)t>n, \mu, \nu \in\mathcal M(\Rn)$ and 
%$\mu(B(x,r))\leq r^s, \nu(B(x,r))\leq r^t$ for $x\in\Rn, r>0$, then
%\begin{equation}\label{eq16b}
%\iint S_{g\#}(\mu\times\nu)(x)^2\,dx\,d\theta_ng \lesssim \mu(\Rn)\nu(\Rn),
%\end{equation}
%with the implicit constant independent of $\mu$ and $\nu$, assuming $A\subset B(0,1)$. The arguments described above give this too.

\section{Level sets and intersections}
The following results with positive lower density assumptions were proven in \cite{M7}. It turned out that those assumptions can be removed quite easily, as explained below.

The estimate \eqref{eq16} can be used to derive information on the Hausdorff dimension of the level sets of $S_g$, and hence, by \eqref{eq16.2}, of intersections. We shall first discuss a more general version of this principle: a quantitative projection theorem leads to estimates of the Hausdorff dimension of level sets. This is also how in \cite[Chapter 10]{M4} the proof for Marstrand's section theorem \ref{mar2} runs.

We consider the following general setting. Let $P_{\lambda}:\Rn\to\R^m, \lambda\in\Lambda,$ be orthogonal projections, where $\Lambda$ is a compact metric space.  Suppose that $\lambda\mapsto P_{\lambda}x$ is  continuous for every $x\in\Rn$. Let also $\omega$ be a finite non-zero Borel measure on $\Lambda$. We denote by $D(\mu,\cdot)$ the Radon-Nikodym derivative of a measure $\mu$ on $\R^m$.

\begin{thm}\label{level}
Let $s>m$. Suppose that there exists a positive number $C$ such that $P_{\lambda\sharp}\mu\ll\mathcal L^m$ for $\omega$ almost all $\lambda\in \Lambda$ and 

\begin{equation}\label{L2}
\iint D(P_{\lambda\sharp}\mu,u)^2\,d\mathcal L^mu\,d\omega\lambda < C
\end{equation}
whenever $\mu\in\mathcal M(B^n(0,1))$ is such that  $\mu(B(x,r))\leq r^s$ for $x\in \Rn, r>0$. 

If $A\subset\R^n$ is $\mathcal H^s$ measurable, $0<\mathcal H^s(A)<\infty$ and $\theta^s_{\ast}(A,x)>0$ (recall \eqref{eq3}) for $\mathcal H^s$ almost all $x\in A$, then 
%for $\mathcal H^s\times\omega$ almost all $(x,\lambda)\in A\times\Lambda$,
%\begin{equation}\label{e17}
%\dim P_{\lambda}^{-1}\{P_{\lambda}x\}\cap A =s-m,
%\end{equation}
for $\omega$ almost all $\lambda\in \Lambda$,
\begin{equation}\label{e18}
\mathcal L^m(\{u\in\R^m: \dim P_{\lambda}^{-1}\{u\}\cap A = s-m\}) > 0.
\end{equation}
\end{thm}

\begin{thm}\label{levela}
Let $s>m$. Suppose that there exists a positive number $C$ such that $P_{\lambda\sharp}\mu\ll\mathcal L^m$ for $\omega$ almost all $\lambda\in \Lambda$ and 

\begin{equation}\label{L2a}
\iint D(P_{\lambda\sharp}\mu,u)^2\,d\mathcal L^mu\,d\omega\lambda < C\mu(B^n(0,1))
\end{equation}
whenever $\mu\in\mathcal M(B^n(0,1))$ is such that  $\mu(B(x,r))\leq r^s$ for $x\in \Rn, r>0$. 

If $A\subset\R^n$ is $\mathcal H^s$ measurable, $0<\mathcal H^s(A)<\infty$, then 
%for $\mathcal H^s\times\omega$ almost all $(x,\lambda)\in A\times\Lambda$,
%\begin{equation}\label{e17}
%\dim P_{\lambda}^{-1}\{P_{\lambda}x\}\cap A =s-m,
%\end{equation}
for $\omega$ almost all $\lambda\in \Lambda$,
\begin{equation}\label{e18a}
\mathcal L^m(\{u\in\R^m: \dim P_{\lambda}^{-1}\{u\}\cap A = s-m\}) > 0.
\end{equation}
\end{thm}

For an application to intersections we shall need the following product set version of Theorem \ref{levela}, Theorem \ref{level1a}. There $P_{\lambda}:\Rn\times\R^p\to\R^m, \lambda\in\Lambda, m<n+p,$ are orthogonal projections with the same assumptions as before.

\begin{thm}\label{level1}
Let $s,t>0$ with $s+t>m$. Suppose that there exists a positive number $C$ such that $P_{\lambda\sharp}(\mu\times\nu)\ll\mathcal L^m$ for $\omega$ almost all $\lambda\in \Lambda$ and 

\begin{equation}\label{L3}
\iint D(P_{\lambda\sharp}(\mu\times\nu),u)^2\,d\mathcal L^mu\,d\omega\lambda < C
\end{equation}
whenever $\mu\in\mathcal M(B^n(0,1)), \nu\in\mathcal M(B^p(0,1))$ are such that $\mu(B(x,r))\leq r^s$ for $x\in \Rn, r>0$, and $\nu(B(y,r))\leq r^t$ for $y\in \R^p, r>0$. 

If $A\subset\R^n$ is $\mathcal H^s$ measurable, $0<\mathcal H^s(A)<\infty$,\ $B\subset\R^p$ is $\mathcal H^t$ measurable, $0<\mathcal H^t(B)<\infty$, $\theta^s_{\ast}(A,x)>0$ for $\mathcal H^s$ almost all $x\in A$, and $\theta^t_{\ast}(B,y)>0$ for $\mathcal H^t$ almost all $y\in B$, then 
%for $\mathcal H^s\times \mathcal H^t\times\omega$ almost all $(x,y,\lambda)\in A\times B\times\Lambda$,
%\begin{equation}\label{pr1}
%\dim P_{\lambda}^{-1}\{P_{\lambda}(x,y)\}\cap (A\times B) =s+t-m,
%\end{equation}
for $\omega$ almost all $\lambda\in \Lambda$,
\begin{equation}\label{pr2}
\mathcal L^m(\{u\in\R^m: \dim P_{\lambda}^{-1}\{u\}\cap (A\times B) = s+t-m\}) > 0.
\end{equation}
\end{thm}

\begin{thm}\label{level1a}
Let $s,t>0$ with $s+t>m$. Suppose that there exists a positive number $C$ such that $P_{\lambda\sharp}(\mu\times\nu)\ll\mathcal L^m$ for $\omega$ almost all $\lambda\in \Lambda$ and 

\begin{equation}\label{L3a}
\iint D(P_{\lambda\sharp}(\mu\times\nu),u)^2\,d\mathcal L^mu\,d\omega\lambda < C\mu(B^n(0,1))\nu(B^p(0,1))
\end{equation}
whenever $\mu\in\mathcal M(B^n(0,1)), \nu\in\mathcal M(B^p(0,1))$ are such that $\mu(B(x,r))\leq r^s$ for $x\in \Rn, r>0$, and $\nu(B(y,r))\leq r^t$ for $y\in \R^p, r>0$. 

If $A\subset\R^n$ is $\mathcal H^s$ measurable, $0<\mathcal H^s(A)<\infty$,\ $B\subset\R^p$ is $\mathcal H^t$ measurable, $0<\mathcal H^t(B)<\infty$, then 
%for $\mathcal H^s\times \mathcal H^t\times\omega$ almost all $(x,y,\lambda)\in A\times B\times\Lambda$,
%\begin{equation}\label{pr1}
%\dim P_{\lambda}^{-1}\{P_{\lambda}(x,y)\}\cap (A\times B) =s+t-m,
%\end{equation}
for $\omega$ almost all $\lambda\in \Lambda$,
\begin{equation}\label{pr2a}
\mathcal L^m(\{u\in\R^m: \dim P_{\lambda}^{-1}\{u\}\cap (A\times B) = s+t-m\}) > 0.
\end{equation}
\end{thm}
I give a few words about the proofs of Theorems \ref{level} and \ref{levela}. First, notice that $D(P_{\lambda\sharp}(\mu),u)$ is given by
$$D(P_{\lambda\sharp}\mu,u) = \lim_{\delta\to 0}\mathcal L^m(B(0,1))^{-1}\delta^{-m} \mu(\{y:|P_{\lambda}(y)-u|\leq\delta\}).$$

Let $\mu$ be the restriction of $\mathcal H^s$ to a subset of $A$ so that $\mu$ satisfies the Frostman $s$ condition. Then \eqref{L2} is applied  to the measures

$$\mu_{a,r,\delta}=r^{-s}T_{a,r\sharp}(\mu_{\delta}\restrict B(a,r))\in\mathcal M(B(0,1)), a\in\Rn, r>0, \delta>0,$$
where $\mu_{\delta}(B)=\delta^{-n}\int_B\mu(B(x,r))\,d\mathcal L^n x$,\ $T_{a,r}(x)=(x-a)/r$ is the blow-up map and $\mu_{\delta}\restrict B(a,r)$ is the restriction of $\mu_{\delta}$ to $B(a,r)$. This leads for almost all $x\in A, \lambda\in\Lambda$, to 
\begin{equation}\label{e8}
\lim_{r\to 0}\liminf_{\delta\to 0}r^{-t}\delta^{-m} \mu(\{y\in B(x,r):|P_{\lambda}(y-x)|\leq\delta\}) = 0,
\end{equation}
which is a Frostman type condition along the level sets of the $P_{\lambda}$. With some further work it leads to \eqref{e18}. The proof of Theorem \ref{level1} is similar.

The proof of Theorem \ref{levela} is the same as that of Theorem \ref{level} observing that the assumption of positive lower density is only used to get the last line of page 394 in \cite{M7}, and this, with $\mu(B_{j,i})$ replaced by $\mu(2B_{j,i})$, now follows from the assumption \eqref{L2a}. Here $\mu(2B_{j,i})$ is fine, since a few lines earlier in \cite{M7} we can require that the balls $2B_{j,i}$ too have bounded overlap.

Theorem \ref{level1a} together with the quantitative version \eqref{eq16b} of Theorem \ref{thm4} and with \eqref{eq16.2} can be applied to the projections $S_g$ to obtain the following result on the Hausdorff dimension of intersections:

\begin{thm}\label{inter}
Let $s,t>0$ with $s+(n-1)t/n>n$ and let $A\subset\R^n$ be $\mathcal H^s$ measurable with $0<\mathcal H^s(A)<\infty$, and let $B\subset\R^n$ be $\mathcal H^t$ measurable with $0<\mathcal H^t(B)<\infty$. Then 
%for $\mathcal H^s\times\mathcal H^t\times\theta_n$ almost all $(x,y,g)\in A\times B\times O(n)$,
%\begin{equation}\label{e5}
%\dim  A\cap (g(B-y)+x) = s+t-n,
%\end{equation}
for $\theta_n$ almost all $g\in O(n)$,
\begin{equation}\label{e6}
\mathcal L^n(\{z\in\R^n:\dim A\cap (g(B)+z) \geq s+t-n\}) > 0.
\end{equation}
\end{thm}

Harris \cite{H2} proved the analogue of the part \eqref{mareq12} of Theorem \ref{mar1} for vertical projections in the Heisenberg group $\H^1$. He applied this with a method partially similar to that of the proof of Theorem \ref{levela} to get the analogue of the intersection theorem \ref{mar2}. This lead me to remove the lower density assumptions.

\section{Exceptional set estimates}\label{Excestimates}

Recall the exceptional set estimates for orthogonal projections and for intersections with planes from Chapter \ref{Projections}. Now we discuss some similar results from \cite{M6} for intersections.

First we have an exceptional set estimate related to Theorem \ref{thm2}. But we need a bit stronger assumption: the sum of the dimensions is required to be  bigger than $n+1$, rather than just one of the sets having dimension bigger than $(n+1)/2$. Recall that the dimension of $O(n)$ is $n(n-1)/2$.

\begin{thm}\label{extheo1}
Let $s$ and $t$ be positive numbers with $s+t > n+1$. Let $A$ and $B$ be Borel subsets of $\Rn$ with $\mathcal H^s(A)>0$ and $\mathcal H^t(B)>0$. Then 
there is $E\subset O(n)$ such that 
$$\dim E\leq n(n-1)/2-(s+t-(n+1))$$ 
and for  $g\in O(n)\setminus E$,
\begin{equation}\label{eq18}
\mathcal L^n(\{z\in\Rn: \dim A\cap (g(B)+z)\geq s+t-n\})>0.
\end{equation}
\end{thm}

The proof is based on the Fourier transform and the convolution approximation mentioned in Section \ref{methods}. Instead of $\theta_n$ one uses a Frostman measure $\theta$ on the exceptional set $E$: if $\a>(n-1)(n-2)/2$ is such that $\theta(B(g,r))\leq r^{\a}$ for all $g\in O(n)$ and $r>0$, then for $x,z\in\Rn\setminus\{0\}, r>0$,
\begin{equation}
\theta(\{g:|x-g(z)|< r\})\lesssim (r/|z|)^{\a-(n-1)(n-2)/2}.
\end{equation}
This replaces the inequality \eqref{theta}.

In the case where one of the sets has small dimension we have the following improvement of Theorem \ref{extheo1}:

\begin{thm}\label{theo6}
Let $A$ and $B$ be Borel subsets of $\Rn$  and suppose that $\dim A\leq (n-1)/2$. If $0<u<\dim A + \dim B - n$, then there is $E\subset O(n)$ with 
$$\dim E\leq n(n-1)/2-u$$ 
such that for  $g\in O(n)\setminus E$,
\begin{equation}\label{eq19}
\mathcal L^n(\{z\in\Rn: \dim A\cap (g(B)+z)\geq u\})>0.
\end{equation}
\end{thm}

The last decay estimate in \eqref{dz} of spherical averages is essential for the proof. The reason why the assumption $\dim A\leq (n-1)/2$ leads to a better result is that that estimate in \eqref{dz} is stronger than the others. For $\dim A > (n-1)/2$ the inequalities \eqref{dz} would only give weaker results with $u$ replaced by a smaller number, see \cite[Section 4]{M6}.

If one of the sets supports a measure with sufficiently fast decay of the averages $\sigma(\mu)(r)$, we can improve the estimate of Theorem \ref{extheo1}. Then the results even hold without any rotations provided the dimensions are big enough. In particular, we have the following result in case one of the sets is a Salem set. By definition, $A$ is a Salem set if for every $0<s<\dim A$ there is $\mu\in\mathcal M(A)$ such that $|\widehat{\mu}(x)|^2\lesssim |x|^{-s}$. A discussion on Salem sets can be found, for example, in \cite{M5}, Section 3.6. 

\begin{thm}\label{theo4}
Let $A$ and $B$ be Borel subsets of $\Rn$  and suppose that $A$ is a Salem set. Suppose that $0<u<\dim A + \dim B - n$.

(a) If $\dim A+\dim B>2n-1$, then 
\begin{equation}\label{eq20}
\mathcal L^n(\{z\in\Rn: \dim A\cap (B+z)\geq u\})>0.
\end{equation}
(b) If $\dim A+\dim B\leq 2n-1$, then there is $E\subset O(n)$ with 
$$\dim E\leq n(n-1)/2-u$$ 
such that for  $g\in O(n)\setminus E$,
\begin{equation}\label{eq21}
\mathcal L^n(\{z\in\Rn: \dim A\cap (g(B)+z)\geq u\})>0.
\end{equation}
\end{thm}

Could this hold for general sets, perhaps in the form that $\dim E = 0$, if $\dim A+\dim B>2n-1$? It follows from Theorem \ref{exmar2} that this is true if one of the sets is a plane. In $\R^2$ a slightly stronger question reads: if $s+t>2$ and  $A$ and $B$ are Borel subsets of $\R^2$ with $\mathcal H^s(A)>0$ and $\mathcal H^t(B)>0$,  is there $E\subset O(2)$ with $\dim E=0$, if $s+t\geq 3$, and 
$\dim E\leq 3-s-t$, if $s+t\leq 3$,  
such that for  $g\in O(2)\setminus E$,
$$\mathcal L^2(\{z\in\R^2: \dim A\cap (g(B)+z)\geq s+t-2\})>0?$$

\section{Some relations to the distance set problem}\label{distset}

There are some connections of this topic to Falconer's distance set problem. For general discussion and references, see for example \cite{M5}. Falconer showed in \cite{F2} that for a Borel set $A\subset\Rn$ the distance set $\{|x-y|: x,y\in A\}$ has positive Lebesgue measure if $\dim A > (n+1)/2$. We had the same condition in Theorem \ref{thm2}. Also for distance sets it is expected that $\dim A>n/2$ should be enough.

When $n=2$ Wolff \cite{W} improved $3/2$ to $4/3$ using \eqref{dz}.  Observe that when $\dim A = \dim B$, the assumption $\dim A + \dim B/2 > 2$ in Theorem \ref{inter} becomes $\dim A > 4/3$ and it is the same as Wolff's. This is no coincidence: both results use Wolff's estimate \eqref{dz}.
%In \cite{IL} Iosevich and Liu improved distance set results of the time for product sets with rather simple arguments. 

The proofs of distance set results often involve the distance measure $\delta(\mu)$ of a measure $\mu$ defined by

$$\delta(\mu)(B) = \mu\times\mu(\{(x,y): |x-y|\in B\}),\ B\subset\R.$$

The crucial estimate \eqref{eq14} means that $\delta(\mu)$ is absolutely continuous with bounded density if $I_{(n+1)/2}(\mu)<\infty$. Hence it yields Falconer's result. As mentioned before we cannot hope to get bounded density when $s<(n+1)/2$, at least when $n=2$ or $3$. In many of the later improvements one verifies absolute continuity with $L^2$  density. For example, Wolff showed that $\delta(\mu)\in L^2(\R)$, if $I_{s}(\mu)<\infty$ for some $s>4/3$. To do this he used decay estimates for the spherical averages $\sigma(\mu)(r)$ and proved \eqref{dz} for $n=2$. The proofs of the most recent,  and so far the best known, distance set results in  \cite{DZ}, \cite{GIOW}, \cite{DGOWWZ} and \cite{DIOWZ} are quite involved using deep harmonic analysis techniques; restriction and decoupling. In the plane the result of \cite{GIOW} says that the distance set of $A$ has positive Lebesgue measure if $\dim A > 5/4$. See Shmerkin's survey \cite{S2} for the distance set and related problems.

Distance measures are related to the projections $S_g$ by the following:

\begin{equation}\label{distproj}
\iint D(S_{g\#}(\mu\times\nu))(z)^2\,d\mathcal L^nz\,d\theta_ng = c\int \delta(\mu)(t)\delta(\nu)(t)t^{1-n}\,dt,\end{equation}
at least if $\mu$ and $\nu$ are smooth functions with compact support, see \cite[Section 5.2]{M8}.

Since by an example in \cite{GIOW}, when $n=2$, for any $s<4/3$,\  $I_{s}(\mu)<\infty$ is not enough for $\delta(\mu)$ to be in $L^2$, probably, because of \eqref{distproj}, it is not enough for  $S_{g\#}(\mu\times\mu)$ to be in $L^2$. But  in \cite{GIOW} it was shown that if $I_{s}(\mu)<\infty$ for some $s>5/4$, there is a complex valued modification of $\mu$ with good $L^2$ behaviour. In higher even dimensions similar results were proven in \cite{DIOWZ} with $n/2+1/4$ in place of $5/4$. Could those methods be used to show, for instance, that if $n=2$ and $\dim A = \dim B > 5/4$, then $\mathcal L^2(\mathcal S_g(A\times B))>0$ for almost all $g\in O(2)$?  
%One problem is that for distance sets one can split the measure to two parts with disjoint supports and only consider distances between points in the different supports. So one need not consider arbitrarily small distances, and the authors of \cite{GIOW} and in \cite{DIOWZ} seem to use this essentially. Here such reduction may not be possible.

%\section{Other dimensions}

\vspace{1cm}
\begin{footnotesize}
{\sc Department of Mathematics and Statistics,
P.O. Box 68,  FI-00014 University of Helsinki, Finland}\\
\emph{E-mail address:} 
\verb"pertti.mattila@helsinki.fi" 

\end{footnotesize}

\end{document}